\documentclass[12pt]{article}
\author{Thomas Morrill\footnote{Supported by Australian Research Council Discovery Project DP160100932.}\\ School of Science\\ The University of New South Wales Canberra, Australia \\ t.morrill@adfa.edu.au\and Dave Platt$^*$\\ School of Mathematics, \\ University of Bristol, Bristol, UK\\dave.platt@bris.ac.uk\and  Tim Trudgian$^{*}$\footnote{Supported by Australian Research Council Future Fellowship FT160100094.} \\
School of Science\\ The University of New South Wales Canberra, Australia \\
  t.trudgian@adfa.edu.au}
  \title{Sign changes in the prime number theorem}
\usepackage{indentfirst}
\usepackage{url}
\usepackage{enumerate}
\usepackage{amsthm}
\usepackage{amsmath}
\usepackage{comment}
\usepackage{fullpage}
\usepackage{amssymb}
\usepackage{booktabs}

\newtheorem{thm}{Theorem}

\newtheorem{Lem}{Lemma}

\newcommand{\dif}{\textrm d}
\newcommand{\e}{\textrm e}
\newcommand{\Ei}{\textrm{Ei}}

\begin{document}
\maketitle
\begin{abstract}
\noindent 
Let $V(T)$ denote the number of sign changes in $\psi(x) - x$ for $x\in[1, T]$. We show that $\liminf_{T\rightarrow\infty} V(T)/\log T\geq \gamma_{1}/\pi +  1.867\cdot 10^{-30}$, where $\gamma_{1} = 14.13\ldots$ denotes the ordinate of the lowest-lying non-trivial zero of the Riemann zeta-function. This improves on a long-standing result by Kaczorowski.
\end{abstract}
\section{Introduction and statement of main result}
\noindent
Let $\pi(x)$ denote the number of primes not exceeding $x$. By the prime number theorem, we have $\pi(x) \sim x/\log x$. Four equivalent versions of this are
\begin{equation}\label{rum}
\pi(x) \sim \textrm{li}(x), \quad \Pi(x) \sim \textrm{li}(x), \quad \psi(x) \sim x, \quad \theta(x)\sim x,
\end{equation}
where $\textrm{li}(x)$ denotes the logarithmic integral, and where
$$\Pi(x) = \pi(x) + \frac{1}{2} \pi(x^{1/2}) + \frac{1}{3} \pi(x^{1/3})\ldots, \quad \psi(x) = \sum_{p^{m} \leq x} \log p, \quad \theta(x) = \sum_{p\leq x} \log p.$$
The finer details of the relations in (\ref{rum}) have been the cause of much investigation. Following Kaczorowski \cite{Kacs1}, for $1\leq i \leq 4$ define $\Delta_{i}(x)$ as the difference of the right and left sides in each of the relations in (\ref{rum}). There are infinitely many changes of sign for each of the $\Delta_{i}(x)$: this result goes back to Littlewood \cite{Littlewood}. Two problems naturally arise: what is the first sign change, and how frequently do these sign changes occur?

The first sign change of $\Delta_{2}(x)$ is at $x=31$; the first sign change of $\Delta_{3}(x)$ is at $x=19$. Locating the first sign changes of $\Delta_{1}(x)$ and $\Delta_{4}(x)$ appears a hopeless endeavour. The history of the former is rich --- see \cite{Saouter} for more details; the latter has not attracted much attention --- see \cite{PT} for the best known upper bound on the first sign change.

The second problem --- the frequency with which sign changes occur --- has been addressed by many authors, perhaps most extensively by Kaczorowski in two series of articles beginning \cite{Kack1} and \cite{Kacs1}. These articles contain meticulous references to earlier work: we refer the reader to those for more details. For $1\leq i \leq 4$, let $V_{i}(T)$ denote the number of sign changes of $\Delta_{i}(x)$ with $x\in[1, T]$. Hereafter we shall only be concerned with $T$ sufficiently large.

In \cite{Kacs1} it is shown that $V_{i}(T)\geq \gamma_{1}/(4\pi) \log T$ for $i=2,3$, where $\gamma_{1} = 14.13\ldots$ is the imaginary part of the lowest-lying non-trivial zero of the Riemann zeta-function. Complementary results are proved in \cite{Kacs2}, in which Kaczorowski improves upon a result by P\'{o}lya \cite{Polya} by showing that $V_{1}(T) \geq c_{1}\log T$ for an (ineffective) constant $c_{1}$. He also proves that
\begin{equation}\label{whisky}
\liminf_{T\rightarrow\infty} \frac{V_{3}(T)}{\log T} \geq \frac{\gamma^{*}}{\pi},
\end{equation}
where $\gamma^{*}$ is defined as follows. Let $\Theta$ denote the supremum over $\Re(\rho)$ where $\rho$ ranges over all zeroes of $\zeta(s)$. If there are any zeroes $\rho = \beta + i\gamma$ with $\beta = \Theta$ we shall define $\gamma^{*}$ as the least positive $\gamma$; otherwise $\gamma^{*} =\infty$. By the work of the second author \cite{Platt2017}, if the Riemann hypothesis is false then $\gamma^{*}> 3\cdot 10^{10}$. If the Riemann hypothesis is true, then $\gamma^{*} = \gamma_{1}$. Hence, the result in (\ref{whisky}) holds with $\gamma^{*}$ replaced by $\gamma_{1}$.

In Theorem 1.2 in \cite{Kack5} Kaczorowski makes a small improvement on (\ref{whisky}), showing that
\begin{equation}\label{bourbon}
\liminf_{T\rightarrow\infty} \frac{V_{3}(T)}{\log T} \geq \frac{\gamma_{1}}{\pi} + 10^{-250}.
\end{equation}
It is stated that no serious effort was expended in producing the quantity $10^{-250}$, and that this constant ought to be improvable.   With such a fundamental question, it seems reasonable to see what  theoretical and computational power can be used to replace $10^{-250}$ in (\ref{bourbon}) with a larger constant. Of course, the conjectured rate of growth of $V_{3}(T)$, roughly as large as $\sqrt{T}$ --- see, e.g., \cite{vortex} --- appears hopelessly beyond current methods.

In this paper we follow Kaczorowski's method, making some theoretical and computational improvements, and prove the following theorem.
\begin{thm}\label{papageno}
We have
\begin{equation*}
\liminf_{T\rightarrow\infty}\frac{V_{3}(T)}{\log T} \geq \frac{\gamma_{1}}{\pi} +1.867\cdot 10^{-30}.
\end{equation*}
\end{thm}
We present an outline of Kaczorowski's method, with some minor improvements, in \S \ref{boat}. We detail our computational approach in \S \ref{yacht} and prove Theorem \ref{papageno}.

\section{Kaczorowski's method}\label{boat}
For $\Re z>0$ and $\Im z >0$ define 
\begin{equation} \label{ground}
	F(z) = e^{-z/2} \sum_{\substack{\varrho = \beta + i \gamma \\ \gamma > 0} } \frac{e^{\varrho z}}{\varrho},
\end{equation}
where the sum is taken over nontrivial zeroes of $\zeta(s)$.
Counting the zeroes of this function gives an improved lower bound on \eqref{whisky}.

\begin{thm}[Corollary 1.1 \cite{Kack5}] \label{Kac-cor}
	Assume the Riemann hypothesis.
	Then,
	\begin{align*}
		\liminf_{T\rightarrow\infty}  \frac{V_3(T)}{\log T} \geq \frac{\gamma_1}{\pi} + 2 \varkappa,
	\end{align*}
	where
	\begin{align*}
		\varkappa = \lim_{Y \to 0^+} \lim_{T \to \infty} \frac{1}{T}\ \# \{z = x + iy\; \;\textrm{such that}\;\;  F(z) = 0,\ 0 < x < T,\ y \geq Y \}.
	\end{align*}
\end{thm}

Let
\begin{equation*} \label{roast}
	F_N(z) =  \sum_{n = 1}^{N} \frac{e^{i \gamma_{n} z}}{1/2 +  i \gamma_{n}} =: \sum_{i = 1}^{N}  a_n e^{i w_n z}.
\end{equation*}
where $\{\gamma_n\}_{n\geq 1}$ is the increasing sequence of positive imaginary parts of zeroes of $\zeta(s)$.
For a zero $\xi =  u + iv$ of $F$ and a Jordan curve $C$ encircling $\xi$, define
\begin{equation*}
\begin{split}
	\alpha &= \inf_C F_N(z),
	\quad y_0 = \inf_C \Im z,
	\quad a = \sum_{n = 1}^{N} \left| a_n e^{- w_n y_0} \right|,
	\quad b = \sum_{n = N+1}^{\infty} \left| a_n e^{- w_n y_0} \right|,\\
	x_{0} &= \inf_{C} \Re z, \quad x_{1} = \sup_{C} \Re z.
	\end{split}
\end{equation*}
With these definitions, Kaczorowski is able to apply Dirichlet's theorem on diophantine approximation to obtain a lower bound on $\varkappa$.
\begin{Lem}[Theorem 7.1 in  \cite{Kack5}] \label{Kac-thm}
	If $ \alpha < 3b$ and $x_{1} - x_{0} <1$ then $\varkappa \geq q_0^{-N}$, where
\begin{equation}\label{sarastro}
		q_0 =
		\bigg[
			\frac{4 \pi a}{\alpha - 3b}
		\bigg] + 1.
\end{equation}
\end{Lem}

We improve this slightly. The second displayed equation on page 55 of \cite{Kack5} bounds $F(\xi + n_{l})$ where $\xi$ is a zero of $F$ and $n_{l}$ is an integer. This is approximated by sums up to $N$, in which Dirichlet's theorem is used, and the tail pieces of the sums, in which the terms are estimated trivially. We can use the location of the zero $\xi = u + iw$, and not the edge of the contour, for a more precise bound on these sums. Also, for $|\theta| < 1/q_{0}$ we have
$$|e^{ 2\pi i \theta} - e^{-2 \pi i \theta}| = |2 \sin \pi \theta| \leq 2 \sin \pi/q_{0} \leq 2\pi/q_{0}.$$
 We retain the first inequality (whereas Kaczorowski uses the second) and, using
$$ a_{w} = \sum_{n=1}^{N} |a_{n}| e^{-w_{n} w}, \quad b_{w} = \sum_{n=N+1}^{\infty} |a_{n}| e^{-w_{n} w},$$  we obtain our condition on $q_{0}$ in the following lemma.
\begin{Lem} \label{Kac-thm2}
	If $ \alpha > b + 2 b_{w}$ and $x_{1} - x_{0} <1$ then $\varkappa \geq q_0^{-N}$, where $q_{0}\in\mathbb{Z}_{>0}$ is such that
\begin{equation}\label{pamina}
2a_{w} \sin (\pi/q_{0}) + \frac{2\pi a}{q_{0}} \leq \alpha - b - 2b_{w}.
\end{equation}
\end{Lem}

We shall detail in the next section that we are able to makes gains using our criterion (\ref{pamina}) that were unobtainable using the earlier result (\ref{sarastro}).

\section{Proof of Theorem \ref{papageno}}\label{yacht}

The computation to prove Theorem \ref{papageno} proceeds in two parts. This first is a non-rigorous search for promising zeroes of $F$ together with an indication of the value of $\varkappa$ that will result from each such zero. The second part takes the best candidate from the first step and rigorously computes $\varkappa$ from scratch.

\subsection{Computation}\label{quartet}
To handle sums over non-trivial zeroes of $\zeta$ we will use the following result due to Lehman \cite[Lem.\ 1]{Lehman1960}. 
\begin{Lem}\label{Lehman}
If $\varphi(t)$ is a continuous function which is positive and monotone decreasing for $2\pi\e\leq T_1\leq t\leq T_2$, then for some $\vartheta$ with $|\vartheta|\leq 1$ we have
\begin{equation*}
\sum\limits_{T_1\leq\gamma\leq T_2}\varphi(\gamma)=\frac{1}{2\pi}\int\limits_{T_1}^{T_2}\varphi(t)\log\frac{t}{2\pi} \dif t +\vartheta\left\{4\varphi(T_1)\log T_1+2\int\limits_{T_1}^{T_2}\frac{\varphi(t)}{t} \dif t\right\}.
\end{equation*}
\end{Lem}
We note that, as in \cite{PT2}, we could improve the constants $4$ and $2$ in Lemma \ref{Lehman} using results from \cite{PTZ} and \cite{TrudgianS2}, but the accuracy is already sufficient for our purposes.

We can limit the region of search using the following lemma. A similar strategy was used in the search for zeroes in previous work by the second and third authors \cite[\S 2.2]{PTP}. 
\begin{Lem}\label{F_zeros}
Assume the Riemann hypothesis. Then the function $F(z)$ defined at (\ref{ground}) has no zeroes for $\Im z> 0.084\,1$.
\begin{proof}
Taking $y_0=0.0841$, the contribution from the first zero has absolute value greater than $0.021\,536$. Summing the absolute values of the contributions of the next $999$ zeroes gives less than $0.021\,528$ and by Lemma \ref{Lehman} the rest of the zeroes contribute no more than $5\cdot 10^{-54}$.
\end{proof}
\end{Lem}

We now describe a non-rigorous search for a good candidate zero of $F(z)$ for Lemma~\ref{Kac-thm}. Emprically, we notice that those zeroes with large imaginary part give the best bound and we search for such by applying the Newton--Raphson method to $F_{1000}(z)$. If, during the Newton--Raphson iteration, we find that the imaginary part of $z$ has become negative or has exceeded $0.085$, we can be fairly certain that the algorithm will not converge, and hence we abort. Starting with $z=t+0.04 i$ with the $t$'s spaced $1/10$ apart, we try a maximum of $25$ iterations and then compare the imaginary part to the largest we have seen so far and discard the smaller. Once all the candidate $t$'s have been tried, we then perform another $100$ Newton--Raphson iterations on our best candidate to try to ensure that we have converged.

Having (hopefully) isolated $\xi$, a zero of $F_{1000}$, we shall sweep yet more rigour under the carpet for now and treat it as a zero of $F$. We choose a rectangular contour around $\xi$. We find that the method is fairly insensitive to $x_0$ and $x_1$ so we simply set them to $\Re \xi-0.05$ and $\Re \xi+0.05$ respectively.

We now choose a candidate $N$ (we try $N=10, 11, \ldots$). We then fix $y_1$ by finding the first maximum of $|F_N(\Re \xi +iy)|$ with $y>\Im \xi$. To do this, we solve for $\Im F'_{N}(\Re z + i y)=0$ for $y\in[\Im z,\Im z + 0.1]$ and hope that this corresponds to a zero of $\Re F'_{N}$ as well. Having fixed $y_1$ we fix $y_0$ and $\alpha$ such that $\alpha=|F_{N}(\Re z + i y_0)|=|F_{N}(\Re z + i y_1)|$. Note that here we are assuming that the minimum of $|F_{N}|$ over the contour happens at these two points but we will check this rigorously later. Also we found empirically that there is no advantage to taking a slightly larger $y_0$. Even though this would reduce $a$, $a_w$, $b$ and $b_w$, it does not make up for the decrease in $\alpha$. We can now calculate (still non-rigorously) $a$, $a_w$, $b$ and $b_w$, check that $\alpha>b+2b_w$ and if so, compute the corresponding $q_0$ and $\varkappa$. We continue increasing $N$ while $\varkappa$ increases and take the largest value of $\varkappa$.

We implemented this algorithm in GP/PARI \cite{Batut2000} running at $115$ digits of precision. It takes about $160$ minutes on a single $16$ core node of the University of Bristol's Bluecrystal Phase III \cite{ACRC2015} to search $t\in[0,200\,000]$. We present as Table \ref{zeros_tab} selected results from this process.

  \begin{table}[h]
    \caption{Some promising looking ``zeroes'' of $F_{1000}$}
    \label{zeros_tab}
\begin{center}
  \begin{tabular}{| r | r |}
  \hline
  \multicolumn{1}{|c}{$\Re z$} & 
  \multicolumn{1}{|c|}{$\Im z$}\\ 
  \hline
  $14\,685.51\ldots$ & $0.0798\ldots$\\
  $141\,914.41\ldots$ & $0.0795\ldots$\\
  $52\,206.82\ldots$ & $0.0794\ldots$\\
  $132\,400.21\ldots$ & $0.0787\ldots$\\
  $78\,306.31\ldots$ & $0.0783\ldots$\\
  $153\,566.13\ldots$ & $0.0785\ldots$\\
  \hline
  \end{tabular}
  \end{center}
  \end{table}
\newpage
The proposed parameters for the zero near $\Re z= 14\,685$ were
\begin{itemize}
\item $N=11$
\item $y_0=0.069\,574\,675$
\item $y_1=0.121\,953\,870$
\item $x_0=14\,685.516\,156\,148\,412\,36-0.05$
\item $x_1=14\,685.516\,156\,148\,412\,36+0.05$
\end{itemize}

\subsection{Making the computation rigorous}

We now take the most promising looking ``zero'' of $F_{1000}$ near $14\,685+0.0798i$ and move forward rigorously. This entails
\begin{itemize}
\item Taking account of rounding errors in our floating point computations.
\item Confirming that there is a zero of $F$ where we think there is and obtaining a lower bound for its imaginary part in order to compute $a_w$ and $b_w$.
\item Finding $\alpha$, the minimum of $F_{11}$ on the rectangular contour $C$ with corners at $x_0+y_0i$ and $x_1+y_1i$.
\item Computing $a$, $b$, $a_w$ and $b_w$.
\item Confirming the conditions of Lemma \ref{Kac-thm} are satisfied, and computing $q_0$ and hence $\varkappa$.
\end{itemize}

\subsubsection{Controlling floating point errors} 
To sidestep the need to undertake a potentially painful analysis of the propogation of floating point rounding errors through the calculations we are about to describe, we perform all our computations using ball arithmetic as implemented in ARB \cite{ARB}.

\subsubsection{Confirming the zero of $F$} 
We have (non-rigorously) located a zero of $F_{1000}$. Using the $4\,520$ zeroes of $\zeta$ below $5\,000$ isolated to an absolute precision of $\pm 2^{-102}$ \cite{Platt2017} we perfom one iteration of Newton--Raphson and  confirm that the resulting $\xi_{4520}$ is within $2^{-20}$ of a zero of $F_{4520}$ by a rigorous application of the argument principle around the circular contour centred on $\xi$ with radius $2^{-20}$. We also find that the minimum attained by $|F_{4520}|$ around this contour exceeds $2.9\cdot 10^{-7}$. Given that the Riemann hypothesis holds at least to height $T=3\cdot 10^{10}$ \cite{Platt2017} we have
$$
|F(z)-F_{4520}(z)|\leq\sum\limits_{4520<\gamma<T}\frac{\e^{-\gamma\Im z}}{\gamma}+\e^{\Re z/2}\sum\limits_{\gamma>T}\frac{\e^{-\gamma\Im z}}{\gamma}.
$$
By Lemma \ref{Lehman} with $\varphi(t)=\exp(-0.0798 t)/t$ this is less than $4\cdot 10^{-176}$ so by Rouch\'e's theorem the same circular contour must also contain a zero of $F$. More precisely, there exists a zero of $F$ within 
$$
  \xi=[14\,685.516\,155\,1,\ 14\,685.516\,157\,2]+[0.079\,831\,7,\ 0.079\,833\,8]i.
$$

\subsubsection{Determining $\alpha$}
Again, using ball arithmetic, we divide the rectangular contour $C$ into pieces of length $10^{-6}$, compute an interval for $F_{11}(z)$ on each piece and take the minimum.

\subsubsection{Computing $a$ etc.}
This is now a trivial computation using the first $11$ zeroes to compute $a$ and $a_w$ and the next $4\,509$ to compute $b$ and $b_w$ but adding a (tiny) error for the tail.

\subsubsection{Computing $q_0$}
The results of the above give us
\begin{itemize}
\item $\alpha\geq 0.005\,179\,11$,
\item $a\leq 0.061\,294\,6$,
\item $a_w\leq 0.046\,355\,3$,
\item $b\leq 0.002\,127\,13$, and
\item $b_w\leq 0.000\,895\,455$.
\end{itemize}
We then find that $q_0=536.4$ satisfies Lemma \ref{Kac-thm2} from which we get $2\varkappa\geq 2\cdot 537^{-11}>1.867\cdot 10^{-30}$.

In the course of this work, we have not applied the theory of holomorphic almost periodic functions.
We note that Kaczorowski's argument via diophantine approximation could be replaced by modifying the following theorem of Jessen and Tornehave from \cite{Jessen}.

\begin{thm}[Jessen and Tornehave]
 Let $f(z)$ be a holomorphic almost periodic function. Define
  \begin{align*}
    J_f(y) = \lim_{T \to \infty} \frac{1}{2T} \int_{-T}^T \log |f(x + iy) | \ dx.
  \end{align*}
  If $J_f$ is differentiable at the points $\alpha, \beta$ with $\alpha < \beta$, then the number of zeroes of $f$ in the rectangle
  \begin{align*}
    \{ z  : -T < \Re z < T, \ \alpha < \Im z < \beta\}
  \end{align*}
  is equal to $T/2\pi ( J_f'(\beta) - J_f'(\alpha)) + o(T)$.
\end{thm}

\section*{Acknowledgements}
We should like to acknowledge Jerzy Kaczorowski for his valuable comments, and Jonathan Bober, Ben Green, David Harvey, and Jesse Thorner for spotting a howler in the abstract.

 \end{document}